\documentclass[12pt]{article}

\RequirePackage{amsthm,amsmath,amsfonts,amssymb}
\RequirePackage[colorlinks,citecolor=blue,urlcolor=blue]{hyperref}
\RequirePackage{comment}

\usepackage{tikz}
\usetikzlibrary{arrows.meta}

\theoremstyle{plain}
\newtheorem{theorem}{Theorem}[section]
\newtheorem{lemma}[theorem]{Lemma}
\newtheorem{proposition}[theorem]{Proposition}

\theoremstyle{remark}

\newtheorem{remark}[theorem]{Remark}

\newcommand{\Loc}{\mathrm{Loc}}
\newcommand{\Sym}{\mathrm{Sym}}
\newcommand{\Canon}{\mathrm{Canon}}

\newcommand{\vc}[1]{\bar #1}

\date{}

\begin{document}

\title{The classification of orthogonal arrays {OA}(2048,14,2,7) and some completely regular codes%
\thanks{%
This is the accepted version of the article published in Discrete Math. 347(5) 2024, 113923(1-8),
\url{https://doi.org/10.1016/j.disc.2024.113923}
}
}


\author{Denis~S.~Krotov%
\thanks{Sobolev Institute of Mathematics, Novosibirsk 630090, Russia. e-mail: dk@ieee.org}
}

\maketitle
\newcommand\GeneratorsTable{
\begin{tabular}{l@{\ }l}
\begin{tikzpicture}[
nd/.style={outer sep=0, inner sep=2,right},
arr/.style={{Stealth[length=1.2mm]}-{Stealth[length=1.2mm]}},
]
\draw
 +(0     ,0) node[nd]     {[}
++(1.1ex, 0) node[nd] (0) {0}
++(1.2ex ,0) node[nd] (1) {0}
++(1.2ex ,0) node[nd] (2) {0}
++(1.2ex ,0) node[nd] (3) {0}
++(1.4ex ,0) node[nd] (4) {0}
++(1.2ex ,0) node[nd] (5) {0}
++(1.2ex ,0) node[nd] (6) {0}
++(1.2ex ,0) node[nd] (7) {0}
++(1.5ex ,0) node[nd] (8) {1}
++(1.2ex ,0) node[nd] (9) {1}
++(1.2ex ,0) node[nd] (a) {1}
++(1.2ex ,0) node[nd] (b) {1}
++(1.4ex ,0) node[nd] (c) {1}
++(1.2ex ,0) node[nd] (d) {1}
++(1.2ex ,0) node[nd] (e) {1}
++(1.2ex ,0) node[nd] (f) {1}
 +(1.2ex ,-0.4em) node[nd]{,}
 +(2.2ex ,-0.07em) node[nd, right] {\ $\mathrm{Id}$\,],};
\end{tikzpicture}
&
\begin{tikzpicture}[
nd/.style={outer sep=0, inner sep=2,right},
arr/.style={{Stealth[length=1.2mm]}-{Stealth[length=1.2mm]}},
]
\draw
 +(0     ,0) node[nd]     {[}
++(1.1ex, 0) node[nd] (0) {1}
++(1.2ex ,0) node[nd] (1) {1}
++(1.2ex ,0) node[nd] (2) {1}
++(1.2ex ,0) node[nd] (3) {1}
++(1.4ex ,0) node[nd] (4) {1}
++(1.2ex ,0) node[nd] (5) {1}
++(1.2ex ,0) node[nd] (6) {1}
++(1.2ex ,0) node[nd] (7) {1}
++(1.5ex ,0) node[nd] (8) {0}
++(1.2ex ,0) node[nd] (9) {0}
++(1.2ex ,0) node[nd] (a) {0}
++(1.2ex ,0) node[nd] (b) {0}
++(1.4ex ,0) node[nd] (c) {0}
++(1.2ex ,0) node[nd] (d) {0}
++(1.2ex ,0) node[nd] (e) {0}
++(1.2ex ,0) node[nd] (f) {0}
 +(1.2ex ,-0.4em) node[nd]{,}
 +(2.2ex ,-0.07em) node[nd, right] {\ $\mathrm{Id}$\,],};
\end{tikzpicture}
\\
\begin{tikzpicture}[
nd/.style={outer sep=0, inner sep=2,right},
arr/.style={{Stealth[length=1.2mm]}-{Stealth[length=1.2mm]}},
]
\draw
 +(0     ,0) node[nd]     {[}
++(1.1ex, 0) node[nd] (0) {0}
++(1.2ex ,0) node[nd] (1) {0}
++(1.2ex ,0) node[nd] (2) {0}
++(1.2ex ,0) node[nd] (3) {0}
++(1.4ex ,0) node[nd] (4) {1}
++(1.2ex ,0) node[nd] (5) {1}
++(1.2ex ,0) node[nd] (6) {1}
++(1.2ex ,0) node[nd] (7) {1}
++(1.5ex ,0) node[nd] (8) {0}
++(1.2ex ,0) node[nd] (9) {0}
++(1.2ex ,0) node[nd] (a) {0}
++(1.2ex ,0) node[nd] (b) {0}
++(1.4ex ,0) node[nd] (c) {0}
++(1.2ex ,0) node[nd] (d) {0}
++(1.2ex ,0) node[nd] (e) {1}
++(1.2ex ,0) node[nd] (f) {1}
 +(1.2ex ,-0.4em) node[nd]{,}
 +(2.2ex ,-0.07em) node[nd, right] {(89)(ab)(cd)(ef)\,],};
\draw[arr] (8.south) .. controls +(-70:0.5em) and +(-110:0.5em) .. (9.south);
\draw[arr] (a.south) .. controls +(-70:0.5em) and +(-110:0.5em) .. (b.south);
\draw[arr] (c.south) .. controls +(-70:0.5em) and +(-110:0.5em) .. (d.south);
\draw[arr] (e.south) .. controls +(-70:0.5em) and +(-110:0.5em) .. (f.south);
\end{tikzpicture}
&
\begin{tikzpicture}[
nd/.style={outer sep=0, inner sep=2,right},
arr/.style={{Stealth[length=1.2mm]}-{Stealth[length=1.2mm]}},
]
\draw
 +(0     ,0) node[nd]     {[}
++(1.1ex, 0) node[nd]   (8) {0}
++(1.2ex ,0) node[nd]   (9) {0}
++(1.2ex ,0) node[nd]   (a) {0}
++(1.2ex ,0) node[nd]   (b) {0}
++(1.4ex ,0) node[nd]   (c) {0}
++(1.2ex ,0) node[nd]   (d) {0}
++(1.2ex ,0) node[nd]   (e) {1}
++(1.2ex ,0) node[nd]   (f) {1}
++(1.5ex ,0) node[nd]   (0) {0}
++(1.2ex ,0) node[nd]   (1) {0}
++(1.2ex ,0) node[nd]   (2) {0}
++(1.2ex ,0) node[nd]   (3) {0}
++(1.4ex ,0) node[nd]   (4) {1}
++(1.2ex ,0) node[nd]   (5) {1}
++(1.2ex ,0) node[nd]   (6) {1}
++(1.2ex ,0) node[nd]   (7) {1}
 +(1.2ex ,-0.4em) node[nd]{,}
 +(2.2ex ,-0.07em) node[nd, right] {(01)(23)(45)(67)\,],};
\draw[arr] (8.south) .. controls +(-70:0.5em) and +(-110:0.5em) .. (9.south);
\draw[arr] (a.south) .. controls +(-70:0.5em) and +(-110:0.5em) .. (b.south);
\draw[arr] (c.south) .. controls +(-70:0.5em) and +(-110:0.5em) .. (d.south);
\draw[arr] (e.south) .. controls +(-70:0.5em) and +(-110:0.5em) .. (f.south);
\end{tikzpicture}
 \\[-0.2em]
\begin{tikzpicture}[
nd/.style={outer sep=0, inner sep=2,right},
arr/.style={{Stealth[length=1.2mm]}-{Stealth[length=1.2mm]}},
]
\draw
 +(0     ,0) node[nd]     {[}
++(1.1ex, 0) node[nd] (0) {1}
++(1.2ex ,0) node[nd] (1) {1}
++(1.2ex ,0) node[nd] (2) {0}
++(1.2ex ,0) node[nd] (3) {0}
++(1.4ex ,0) node[nd] (4) {0}
++(1.2ex ,0) node[nd] (5) {0}
++(1.2ex ,0) node[nd] (6) {0}
++(1.2ex ,0) node[nd] (7) {0}
++(1.5ex ,0) node[nd] (8) {1}
++(1.2ex ,0) node[nd] (9) {1}
++(1.2ex ,0) node[nd] (a) {0}
++(1.2ex ,0) node[nd] (b) {0}
++(1.4ex ,0) node[nd] (c) {0}
++(1.2ex ,0) node[nd] (d) {0}
++(1.2ex ,0) node[nd] (e) {0}
++(1.2ex ,0) node[nd] (f) {0}
 +(1.2ex ,-0.4em) node[nd]{,}
 +(2.2ex ,-0.07em) node[nd, right] {(45)(67)(cd)(ef)\,],};
\draw[arr] (4.south) .. controls +(-70:0.5em) and +(-110:0.5em) .. (5.south);
\draw[arr] (6.south) .. controls +(-70:0.5em) and +(-110:0.5em) .. (7.south);
\draw[arr] (c.south) .. controls +(-70:0.5em) and +(-110:0.5em) .. (d.south);
\draw[arr] (e.south) .. controls +(-70:0.5em) and +(-110:0.5em) .. (f.south);
\end{tikzpicture}
 \\[-0.2em]
\begin{tikzpicture}[
nd/.style={outer sep=0, inner sep=2,right},
arr/.style={{Stealth[length=1.2mm]}-{Stealth[length=1.2mm]}},
]
\draw
 +(0     ,0) node[nd]     {[}
++(1.1ex, 0) node[nd] (0) {0}
++(1.2ex ,0) node[nd] (1) {0}
++(1.2ex ,0) node[nd] (2) {1}
++(1.2ex ,0) node[nd] (3) {1}
++(1.4ex ,0) node[nd] (4) {0}
++(1.2ex ,0) node[nd] (5) {0}
++(1.2ex ,0) node[nd] (6) {1}
++(1.2ex ,0) node[nd] (7) {1}
++(1.5ex ,0) node[nd] (8) {0}
++(1.2ex ,0) node[nd] (9) {0}
++(1.2ex ,0) node[nd] (a) {0}
++(1.2ex ,0) node[nd] (b) {0}
++(1.4ex ,0) node[nd] (c) {0}
++(1.2ex ,0) node[nd] (d) {1}
++(1.2ex ,0) node[nd] (e) {0}
++(1.2ex ,0) node[nd] (f) {1}
 +(1.2ex ,-0.4em) node[nd]{,}
 +(2.2ex ,-0.07em) node[nd, right] {(8a)(9b)(ce)(df)\,],};
\draw[arr] (8.south) .. controls +(-70:0.5em) and +(-110:0.5em) .. (a.south);
\draw[arr] (9.south) .. controls +(-70:0.5em) and +(-110:0.5em) .. (b.south);
\draw[arr] (c.south) .. controls +(-70:0.5em) and +(-110:0.5em) .. (e.south);
\draw[arr] (d.south) .. controls +(-70:0.5em) and +(-110:0.5em) .. (f.south);
\end{tikzpicture}
&
\begin{tikzpicture}[
nd/.style={outer sep=0, inner sep=2,right},
arr/.style={{Stealth[length=1.2mm]}-{Stealth[length=1.2mm]}},
]
\draw
 +(0     ,0) node[nd]     {[}
++(1.1ex, 0) node[nd]  (8) {0}
++(1.2ex ,0) node[nd]  (9) {0}
++(1.2ex ,0) node[nd]  (a) {0}
++(1.2ex ,0) node[nd]  (b) {0}
++(1.4ex ,0) node[nd]  (c) {0}
++(1.2ex ,0) node[nd]  (d) {1}
++(1.2ex ,0) node[nd]  (e) {0}
++(1.2ex ,0) node[nd]  (f) {1}
++(1.5ex ,0) node[nd]  (0) {0}
++(1.2ex ,0) node[nd]  (1) {0}
++(1.2ex ,0) node[nd]  (2) {1}
++(1.2ex ,0) node[nd]  (3) {1}
++(1.4ex ,0) node[nd]  (4) {0}
++(1.2ex ,0) node[nd]  (5) {0}
++(1.2ex ,0) node[nd]  (6) {1}
++(1.2ex ,0) node[nd]  (7) {1}
 +(1.2ex ,-0.4em) node[nd]{,}
 +(2.2ex ,-0.07em) node[nd, right] {(02)(13)(46)(57)\,],};
\draw[arr] (8.south) .. controls +(-70:0.5em) and +(-110:0.5em) .. (a.south);
\draw[arr] (9.south) .. controls +(-70:0.5em) and +(-110:0.5em) .. (b.south);
\draw[arr] (c.south) .. controls +(-70:0.5em) and +(-110:0.5em) .. (e.south);
\draw[arr] (d.south) .. controls +(-70:0.5em) and +(-110:0.5em) .. (f.south);
\end{tikzpicture}
 \\[-0.2em]
\begin{tikzpicture}[
nd/.style={outer sep=0, inner sep=2,right},
arr/.style={{Stealth[length=1.2mm]}-{Stealth[length=1.2mm]}},
]
\draw
 +(0     ,0) node[nd]     {[}
++(1.1ex, 0) node[nd] (0) {0}
++(1.2ex ,0) node[nd] (1) {1}
++(1.2ex ,0) node[nd] (2) {0}
++(1.2ex ,0) node[nd] (3) {1}
++(1.4ex ,0) node[nd] (4) {0}
++(1.2ex ,0) node[nd] (5) {1}
++(1.2ex ,0) node[nd] (6) {0}
++(1.2ex ,0) node[nd] (7) {1}
++(1.5ex ,0) node[nd] (8) {0}
++(1.2ex ,0) node[nd] (9) {0}
++(1.2ex ,0) node[nd] (a) {0}
++(1.2ex ,0) node[nd] (b) {1}
++(1.4ex ,0) node[nd] (c) {0}
++(1.2ex ,0) node[nd] (d) {0}
++(1.2ex ,0) node[nd] (e) {0}
++(1.2ex ,0) node[nd] (f) {1}
 +(1.2ex ,-0.4em) node[nd]{,}
 +(2.2ex ,-0.07em) node[nd, right] {(8c)(9d)(ae)(bf)\,],};
\draw[arr] (8.south) .. controls +(-70:0.6em) and +(-110:0.6em) .. (c.south);
\draw[arr] (9.south) .. controls +(-70:0.6em) and +(-110:0.6em) .. (d.south);
\draw[arr] (a.south) .. controls +(-70:0.6em) and +(-110:0.6em) .. (e.south);
\draw[arr] (b.south) .. controls +(-70:0.6em) and +(-110:0.6em) .. (f.south);
\end{tikzpicture}
&
\begin{tikzpicture}[
nd/.style={outer sep=0, inner sep=2,right},
arr/.style={{Stealth[length=1.2mm]}-{Stealth[length=1.2mm]}},
]
\draw
 +(0     ,0) node[nd]     {[}
++(1.1ex, 0) node[nd] (8) {0}               
++(1.2ex ,0) node[nd] (9) {0}               
++(1.2ex ,0) node[nd] (a) {0}               
++(1.2ex ,0) node[nd] (b) {1}               
++(1.4ex ,0) node[nd] (c) {0}               
++(1.2ex ,0) node[nd] (d) {0}               
++(1.2ex ,0) node[nd] (e) {0}               
++(1.2ex ,0) node[nd] (f) {1}               
++(1.5ex ,0) node[nd] (0) {0}
++(1.2ex ,0) node[nd] (1) {1}
++(1.2ex ,0) node[nd] (2) {0}
++(1.2ex ,0) node[nd] (3) {1}
++(1.4ex ,0) node[nd] (4) {0}
++(1.2ex ,0) node[nd] (5) {1}
++(1.2ex ,0) node[nd] (6) {0}
++(1.2ex ,0) node[nd] (7) {1}
 +(1.2ex ,-0.4em) node[nd]{,}
 +(2.2ex ,-0.07em) node[nd, right] {(04)(15)(26)(37)\,],};
\draw[arr] (8.south) .. controls +(-70:0.6em) and +(-110:0.6em) .. (c.south);
\draw[arr] (9.south) .. controls +(-70:0.6em) and +(-110:0.6em) .. (d.south);
\draw[arr] (a.south) .. controls +(-70:0.6em) and +(-110:0.6em) .. (e.south);
\draw[arr] (b.south) .. controls +(-70:0.6em) and +(-110:0.6em) .. (f.south);
\end{tikzpicture}
 \\[-0.2em]
\begin{tikzpicture}[
nd/.style={outer sep=0, inner sep=2,right},
arr/.style={{Stealth[length=1.2mm]}-{Stealth[length=1.2mm]}},
]
\draw
 +(0     ,0) node[nd]     {[}
++(1.1ex, 0) node[nd] (0) {0}
++(1.2ex ,0) node[nd] (1) {0}
++(1.2ex ,0) node[nd] (2) {0}
++(1.2ex ,0) node[nd] (3) {0}
++(1.4ex ,0) node[nd] (4) {1}
++(1.2ex ,0) node[nd] (5) {1}
++(1.2ex ,0) node[nd] (6) {0}
++(1.2ex ,0) node[nd] (7) {0}
++(1.5ex ,0) node[nd] (8) {1}
++(1.2ex ,0) node[nd] (9) {0}
++(1.2ex ,0) node[nd] (a) {1}
++(1.2ex ,0) node[nd] (b) {0}
++(1.4ex ,0) node[nd] (c) {0}
++(1.2ex ,0) node[nd] (d) {0}
++(1.2ex ,0) node[nd] (e) {0}
++(1.2ex ,0) node[nd] (f) {0}
 +(1.2ex ,-0.4em) node[nd]{,}
 +(2.2ex ,-0.07em) node[nd, right] {(23)(67)(ce)(df)\,],};
\draw[arr] (2.south) .. controls +(-70:0.5em) and +(-110:0.5em) .. (3.south);
\draw[arr] (6.south) .. controls +(-70:0.5em) and +(-110:0.5em) .. (7.south);
\draw[arr] (c.south) .. controls +(-70:0.5em) and +(-110:0.5em) .. (e.south);
\draw[arr] (d.south) .. controls +(-70:0.5em) and +(-110:0.5em) .. (f.south);
\end{tikzpicture}
&
\begin{tikzpicture}[
nd/.style={outer sep=0, inner sep=2,right},
arr/.style={{Stealth[length=1.2mm]}-{Stealth[length=1.2mm]}},
]
\draw
 +(0     ,0) node[nd]     {[}
++(1.1ex, 0) node[nd] (8) {1}            
++(1.2ex ,0) node[nd] (9) {0}            
++(1.2ex ,0) node[nd] (a) {1}            
++(1.2ex ,0) node[nd] (b) {0}            
++(1.4ex ,0) node[nd] (c) {0}            
++(1.2ex ,0) node[nd] (d) {0}            
++(1.2ex ,0) node[nd] (e) {0}            
++(1.2ex ,0) node[nd] (f) {0}            
++(1.5ex ,0) node[nd] (0) {0} 
++(1.2ex ,0) node[nd] (1) {0} 
++(1.2ex ,0) node[nd] (2) {0} 
++(1.2ex ,0) node[nd] (3) {0} 
++(1.4ex ,0) node[nd] (4) {1} 
++(1.2ex ,0) node[nd] (5) {1} 
++(1.2ex ,0) node[nd] (6) {0} 
++(1.2ex ,0) node[nd] (7) {0} 
 +(1.2ex ,-0.4em) node[nd]{,}
 +(2.2ex ,-0.07em) node[nd, right] {(46)(57)(ab)(ef)\,].};
\draw[arr] (2.south) .. controls +(-70:0.5em) and +(-110:0.5em) .. (3.south);
\draw[arr] (6.south) .. controls +(-70:0.5em) and +(-110:0.5em) .. (7.south);
\draw[arr] (c.south) .. controls +(-70:0.5em) and +(-110:0.5em) .. (e.south);
\draw[arr] (d.south) .. controls +(-70:0.5em) and +(-110:0.5em) .. (f.south);
\end{tikzpicture}
\end{tabular}
}

\begin{abstract}
We describe the classification of orthogonal arrays OA$(2048,14,2,7)$, or, equivalently, completely regular $\{14;2\}$-codes in the $14$-cube ($30848$ equivalence classes). In particular, we find that there is exactly one almost-OA$(2048,14,2,7{+}1)$, up to equivalence. As derived objects, \linebreak[4] OA$(1024,13,2,6)$ ($202917$ classes) and completely regular $\{12,2;2,12\}$- and $\{14, 12, 2; 2, 12, 14\}$-codes in the $13$- and $14$-cubes, respectively, are also classified.

MSC:  
05B15, 94B25.

Keywords: 
binary orthogonal array, completely regular code, binary 1-perfect code.
\end{abstract}


\section{Introduction}
{O}{rthogonal}
 arrays are combinatorial structures
important both for practical applications
like design of experiments or software testing 
and theoretically, because of many relations with coding theory,
cryptography, design theory, etc., see e.g.~\cite{HSS:OA}.
The classification of orthogonal arrays with given parameters 
is a problem that attracts attention of many researchers,
see e.g.~\cite{BMS:2017:few}, \cite{BulRy:2018}, \cite{SEN:2010:OA}, \cite{PWLL:2021} and the bibliography there.
The main result of the current work is the classification 
of orthogonal arrays OA$(2048,14,2,7)$,
whose cardinality is the largest over all orthogonal arrays
that have ever been computationally classified.

Let us introduce some concepts and notations.
The $n$-cube $Q_n$ is a graph whose vertex set is the set $\{0,1\}^n$ of binary words of length~$n$, which forms a vector space over~$\mathbb{F}_2$;
two such words are adjacent in~$Q_n$ if and only if they differ in exactly one position.
A~binary \emph{orthogonal array} OA$(N,n,2,t)$
is a multiset 
of vertices of~$Q_n$
of cardinality~$N$
such that every subgraph 
isomorphic to~$Q_{n-t}$ 
contains 
$N/2^t$ words from the array.
The considered parameters OA$(2048,14,2,7)$ attain the Friedman bound~\cite{Friedman:92}
\begin{equation}\label{eq:Fried}
N\ge 2^n\Big(1 - \frac{n}{2(t+1)}\Big)
\end{equation}
and hence any such array is simple (without multiple elements) and, moreover, completely regular.
Binary orthogonal arrays with small parameters lying on 
the Friedman bound have been classified in a sequence of works,
see Table~\ref{t:1} (only integer parameters 
with $t<n-1$ and satisfying the necessary condition $t\le 2n/3-1$ \cite{FDF:CorrImmBound} are included).
The classification of OA$(1024,12,2,7)$, OA$(1536,13,2,7)$, 
and OA$(2048,15,2,7)$ was essentially based on such arrays 
being completely regular codes
(OA$(2048,15,2,7)$ are $1$-perfect codes).
In this paper, we further develop the technique
for classification of completely regular codes
to characterize all orthogonal arrays OA$(2048,14,2,7)$.
In particular, we firstly use an ILP solver to speed up the classification of intermediate objects, called local codes,
and partition the class of all OA$(2048,14,2,7)$ 
into two subclasses, depending on containing a special subconfiguration, specific for the considered parameters.
As derived arrays, we also classify  OA$(1024,13,2,6)$.

\begin{table}[ht]
\begin{tabular}{lll||lll}
Parameters &  \multicolumn{2}{l||}{Number of arrays} &
Parameters &  \multicolumn{2}{l}{Number of arrays}
\\ \hline
 OA$(2,3,2,1)$ &  $1$                             &  &    OA$(1,2,2,0)$ &        $1$ &  \\
 OA$(16,6,2,3)$ &  $1$                            &  &    OA$(8,5,2,2)$ &        $1$ &  \\
 OA$(16,7,2,3)$ &  $1$                            &  &    OA$(8,6,2,2)$ &        $1$ &  \\
 OA$(128,9,2,5)$ &  $2$  & \cite{Kirienko2002}    &    OA$(64,8,2,4)$ &       $3$ &   \cite{SEN:2010:OA} \\
 OA$(1024,12,2,7)$ & $16$ &     \cite{KroVor:2020}   &    OA$(512,11,2,6)$ &     $37$ &   \cite{KroVor:2020} \\
 OA$(1536,13,2,7)$ & $1$ &      \cite{Kro:OA13}       &    OA$(768,12,2,6)$ &     $3$ &   \cite{Kro:OA13} \\
 OA$(2048,14,2,7)$ & \underline{$30848$}  & Th.\,\ref{th:main14}                       &    OA$(1024,13,2,6)$ &    \underline{$202917$} &   Th.\,\ref{th:main13} \\
 OA$(2048,15,2,7)$ &  $5983$  &   \cite{OstPot:15}       &    OA$(1024,14,2,6)$ &    $38408$ &  \cite{OstPot:15} \\ 
 OA$(8192,15,2,9)$ &  ?                           &  &     OA$(4096,14,2,8)$ &   ? &  \\ \hline 
\end{tabular}
\caption{Small parameters OA$(N,n,2,t)$ and OA$(N/2,n-1,2,t-1)$,
$N= 2^n(1-n/2(t+1))$, $t \le 2n/3-1$ }
\label{t:1}
\end{table}

A set~$C$ of vertices (code) of 
a graph~$G$
is called
\emph{completely regular}
with \emph{covering radius}~$\rho$
and \emph{intersection array}
$\{b_0,b_1,...,b_{\rho-1};
c_1,...,c_{\rho}\}$
if the partition
$(C=C^{(0)},C^{(1)},...,C^{(\rho)})$
of the vertex set with respect to the distance from~$C$
satisfies the following property:
every vertex from~$C^{(i)}$ 
has exactly~$b_i$ neighbors
in~$C^{(i+1)}$
and~$c_i$ neighbors
in~$C^{(i-1)}$ 
(it is implied that $b_\rho=c_0=0$%
), 
see, e.g., the survey~\cite{BRZ:CR}. A~completely regular code with intersection array 
$\{n;1\}$ in~$Q_n$ is called $1$-perfect.
For a code $C$ in~$Q_n$, the following two operations,
which result in a code in~$Q_{n-1}$, are defined.
The \emph{puncturing} is removing the symbol
in some fixed position from all codewords;
the \emph{shortening} is the same,
but we keep only the codewords 
that contained some fixed symbol, $0$ or $1$, 
in the deleted position.
The \emph{weight}
of a binary word 
is the number of $1$s in it.
By $\mathrm{even}(C)$ and $\mathrm{odd}(C)$,
we denote the even- and odd-weight subcodes of the code~$C$, respectively. The all-zero (all-one) word is denoted by~$\bar0$~($\bar 1$).

Orthogonal arrays attaining bound~\eqref{eq:Fried}
are simple, independent (there are no two adjacent vertices in the array), 
and completely regular with intersection array $A=\{n;2(t+1)-n\}$,
see e.g.~\cite{Kro:Redundancy23}. For brevity,
we will call them $A$-codes 
(in our case, $\{14;2\}$-codes). 
Two codes $C$, $C'$ in~$Q_n$ are 
\emph{equivalent}
if $C=\sigma (C')$ for some automorphism 
of~$Q_n$.

In Section~\ref{s:comp},
we describe the classification technique;
the results are shown in Section~\ref{s:res}.
In Section~\ref{s:prop},  we consider some properties
of the classified OA$(2048,14,2,7)$,
including the property of being 
``the best possible'', almost-OA$(2048,14,2,7+1)$.
The resulting data is available at~\cite{Perfect-related}.

\section{Computation}\label{s:comp}
The classification approach is based on the concept
of local codes.
We say that a set 
$P \subset \{0,1\}^{14}$ 
is an \emph{$r$-local code}
if

\begin{itemize}
  \item[(I)] (\emph{locality condition})
  $P$ consists of words of weight $\le r$;
  \item[(II)] (\emph{up-to-equivalence condition}) $\bar 0$ is not in~$P$;
  \item[(III)] (\emph{exact cover condition}) the neighborhood of every vertex~$\vc v$ of weight less than~$r$
satisfies 
the local condition from 
the definition of 
a $\{14;2\}$-code:
if $\vc v\in P$, 
then $\vc v$ has no neighbors in~$P$;
if $\vc v\not\in P$, then
$\vc v$ has exactly~$2$ neighbors in~$P$;
\item[(IV)] (\emph{boundary inequality condition})
 each word in~$\{0,1\}^{14}$ has at most~$2$ neighbors in~$P$.
\end{itemize}
It is easy to see that up to equivalence
(this motivates the name of condition~(II))
every $\{14;2\}$-code~$C$ includes an $r$-local code~$R$,
defined as the set
of codewords of~$C$
of weight at most~$r$
(i.e., (I) immediately holds).
Indeed, (III) holds 
by the definition
of completely regular code
and because all 
$C$-neighbors of a
vertex~$\vc v$ of weight
$<r$
are in~$R$.
We cannot say the same
for a vertex~$v$ of
``boundary''
weight~$r$ or $r+1$; 
however, our
quotient matrix guarantees
that such a vertex 
cannot have more than~$2$
code neighbors,
which yields~(IV)
(for different quotient matrix, 
one can find different 
boundary inequalities
to reduce the number 
of local codes, see e.g.~\cite{Kro:21066}).
Finally if (II) 
is not satisfied by~$C$
itself, then 
it is satisfied for 
some code 
equivalent to~$C$.

Straightforwardly from definitions,
we have the following key fact.
\begin{lemma}\label{l:lo}
Any $\{14;2\}$-code
that does not contain~$\vc 0$
is a $14$-local code.
In particular, 
any $\{14;2\}$-code
is equivalent to a $14$-local code.
If $C$ is an $r'$-local code
and $R$ consists of all codewords of~$C$
of weight at most~$r$, where $0\le r < r'$,
then $R$ is an $r$-local code.
\end{lemma}

In the last case, we say that
$C$ is an \emph{$r'$-continuation}
of~$R$.

If for each $r$-local code
we can construct all its
$(r+1)$-continuations,
then all local codes and finally 
$\{14;2\}$-codes can be classified
recursively.
However, it turns out that the step
from $2$-local to $3$-local codes
does not finish in a reasonable 
computational time. 
To handle this problem, we
add an intermediate step
based on $(2,3)$-codes,
defined as follows.
A set 
$P \subset \{0,1\}^{14}$ 
is an 
\emph{$(r_0,r_1)$-local code} if
\begin{itemize}
  \item[(0)] 
  $r_1 \in \{r_0,r_0+1\}$ 
  (the case $r_1=r_0+2$
  also has sense but is not used in the current study);
  \item[(I')] (\emph{locality condition})
  each codeword of~$P$
  starts with~$0$ and has weight at most~$r_0$ or starts with~$1$ and has weight at most~$r_1$;
  \item[(II')] (\emph{up-to-equivalence condition}) $P$ contains $10...0$ (and hence does not contain~$\vc 0$);
  \item[(III')] (\emph{exact cover condition}) the neighborhood of every vertex~$\vc v=(v_1,...,v_{14})$ of weight less than~$r_{v_1}$
satisfies 
the local condition from 
the definition of 
a $\{14;2\}$-code:
if $\vc v\in P$, 
then $\vc v$ has no neighbors in~$P$;
if $\vc v\not\in P$, then
$\vc v$ has exactly~$2$ neighbors in~$P$;
\item[(IV)] (\emph{boundary inequality condition}) each word in~$\{0,1\}^{14}$ has at most~$2$ neighbors in~$P$.
\end{itemize}

By a \emph{local code}, we will mean 
an $r$-local code for some~$r$ or 
an $(r_0,r_1)$-local code for some~$(r_0,r_1)$.
Since the first coordinate
plays a special role in the definition of 
$(r_0,r_1)$-local codes,
we define equivalence
of them differently.
Two $(r_0,r_1)$-local codes
$P$ and $P'$ are 
equivalent if
there is a permutation~$\pi$
of coordinates
that fixes the first
coordinate (and hence the codeword $10...0$)
and sends $P$ to $P'$.
In particular, two
equivalent $3$-local codes
can be nonequivalent
as $(3,3)$-local codes.

Now, we describe
the general line of the classification. 
\begin{itemize}
    \item We start with representatives of the $5$ equivalence classes of $(2,2)$-local codes,
    see Section~\ref{s:2loc}.
    \item For each of the $5$ representatives from the step above,
     construct all $(2,3)$-continuations,
     see Section~\ref{s:ex-cov}.
    \item Classify the found $(2,3)$-continuations
    up to equivalence, keep representatives of each
    equivalence class,
     see Section~\ref{s:isom}.
    \item Validate the computations with the orbit-stabilizer theorem (see Section~\ref{s:validation}).
    \item For each of the representatives 
    from the step above,
    construct all $(3,3)$-continuations; 
    classify them up to equivalence; 
    validate.
    \item Classify the found $(3,3)$-local codes up to equivalence as $3$-local codes
    (see Section~\ref{s:33to3}).
    \item For each of the representatives 
    of $(i-1)$-local codes,
    $i=4$, $5$, \ldots, $14$,
    from the step above,
    construct all $i$-continuations; 
    classify them up to equivalence; 
    validate.    
\end{itemize}
Finally we find all nonequivalent 
$\{14;2\}$-codes.

\begin{remark}
   Actually, using the self-complementary
property of a completely regular
code, one can reconstruct a
$\{14;2\}$-code from 
a~$7$-local partition
and even from a~$6$-local partition
(the reconstruction is unique if exists,
but the empiric fact that a $6$-local
partition always has a continuation
is not yet theoretically explained).
However, because of the uniqueness 
of the continuation at last steps, 
this observation does not lead 
to an essential time gain.
\end{remark}

The approach described above
is still time consuming
(but definitely doable in hundreds
core-years),
and to make the search
(approximately four--five times)
faster,
we divide the class of 
$\{14;2\}$-codes
into two subclasses,
called square codes and 
square-free codes,
see Section~\ref{s:sq-sqf}.
With this modified approach,
we only need to continue a small
part of all $3$-local codes,
which essentially quickens 
the most time-consuming step
of the computation.
In several subsections below, 
we describe details of each step 
of the classification.

\subsection[(2,2)-local codes]{$(2,2)$-Local codes}\label{s:2loc}

There is one class of $1$-local codes,
as well as $(1,2)$-local codes,
with the representative 
$\{10...0, 0...01\}$.

\begin{proposition}\label{p:2loc}
There are five classes 
of $(2,2)$-local codes 
(as well as $2$-local codes).
\end{proposition}
\begin{proof}
Consider a $(2,2)$-local code~$P$.
    Apart from the codeword $10...0$,
    there is one another
    weight-$1$ codeword, say $0...01$.
    These two words have no 
    codeword neighbors.
    Each of the $12$ other 
    weight-$1$ words
    has exactly~$2$
    weight-$2$ codeword
    neighbors.
    On the other hand,
    every weight-$2$ codeword
    has exactly~$2$
    weight-$1$ neighbors.
    Considering the $12$
    weight-$1$ non-code words as vertices and the  
    weight-$2$ codewords as edges,
    we have a $2$-regular graph 
    of order~$12$. 
    From~(IV),
    we see that it does not
    have $3$-cycles.
    There are exactly~$5$
    equivalence classes of $2$-regular triangle-free graphs of order~$12$, with cycle
    structures $4+4+4$, $4+8$,
    $5+7$, $6+6$, $12$.
\end{proof}
We call the corresponding $(2,2)$-local codes $L_{4,4,4}$, $L_{4,8}$, $L_{5,7}$, $L_{6,6}$, and~$L_{12}$.

\subsection{Continuation by exact covering}
\label{s:ex-cov}

We first describe how to construct all 
$r$-local codes
$P_{r-1} \cup R_r$ that continue
 a given
$(r-1)$-local code~$P_{r-1}$.
To build such a continuation, we need to find
a set~$R_r$ of vertices of weight~$r$
such that it satisfies~(IV) and the union 
$P_{r-1} \cup R_r$ satisfies~(III)
for all vertices $\vc v$ of weight~$r-1$
(indeed, for all vertices $\vc v$ of weight less than~$r-1$, 
(III) is satisfied because $P_{r-1}$ 
is an $(r-1)$-local code). 
In particular, (III) means that the vertices 
from~$R_r$ are not adjacent to elements
of~$P_{r-1}$. We start with defining the set~$S$ of 
\emph{candidates}, which is the set of all vertices of weight~$r$ that are not adjacent to 
any element of $P_{r-1}$.
We also define the set~$T$ of all weight-$(r-1)$ 
words that are not in~$P_{r-1}$.
Next, we construct a $|T|\times |S|$ $\{0,1\}$-matrix $M$ whose rows and columns 
are indexed by elements
of~$T$ and~$S$ respectively and such that $1$ 
in the entry $(\vc x,\vc y)$ means that $\vc x$ and~$\vc y$ are adjacent. 
For each vertex~$\vc x$ from~$T$, we define 
the \emph{deficiency}~$\delta_{\vc x}$
as $2$ minus the number of its neighbors 
in~$P_{r-1}$. Now, $P_{r-1} \cup R_r$
satisfying~(III) is equivalent to 
each $\vc x$ from~$T$ 
having exactly~$\delta_{\vc x}$
neighbors in~$R_r\subset S$.
Equivalently,
the sum of columns of~$M$ indexed by~$R_r$
equals the deficiency vector $(\delta_{\vc x})_{\vc x\in T}$.
Finding such set of columns is as instance of 
the exact covering problem (which motivates the name of condition~(III)), with multiplicities $\delta_{\vc x}$ of the covered elements.
All solutions can be found by 
Donald Knuth's Algorithm~X with dancing links~\cite{Knuth:DLX},
which is implemented in the library 
\texttt{libexact}~\cite{KasPot08}.
Using this software, 
we can find all solutions
of the exact covering problem,
and it remains to choose those that satisfy~(IV).
Alternatively,
it is possible to modify Algorithm~X
to take into account the additional inequalities 
that come from condition~(IV).
Empirically, this could give some speed gain,
but it was not critical because 
the most of computing time was spent
for filtering the $3$-local codes 
by an ILP solver, as described in the rest of this section.

Now assume that we want to construct all 
$(r,r+1)$-local codes
$P \cup R$ that continue
 a given
$(r,r)$-local code~$P$ 
(constructing $(r+1,r+1)$-continuations
of a $(r,r+1)$-local code is similar).
We need to find
a set~$R$ of vertices of weight~$r+1$ with the first symbol~$1$
such that $P \cup R$ satisfies~(IV) and~(III').
Condition~(III') 
is sufficient to be checked
for all vertices $\vc v$ of weight~$r$ with the first symbol~$1$
(indeed, for all other vertices, 
(III') is satisfied because $P$ 
is an $(r,r)$-local code). 
In particular, (III') means that the vertices 
from~$R$ are not adjacent to elements
of~$P$. We  define the set~$S$ of 
\emph{candidates}, which is the set of 
all vertices of weight~$r+1$ with the first symbol~$1$ 
that are not adjacent to 
any element of~$P$.
We also define the set~$T$ 
of all weight-$r$ 
words that are not in~$P$ and have the first symbol~$1$.
The rest is similar to the case of $r$-local partitions:
we construct a $|T|\times |S|$ matrix, 
solve the corresponding exact covering problem,
and choose only solutions that satisfy~(IV).
One reasonable improvement follows from the possibility to check~(IV) partially
at the step of choosing candidates. 
Indeed, it is possible that some candidate (recall that it has~$1$ in the first position) is adjacent to a vertex with~$0$
in the first position that already has $2$ neighbors in~$P$.
Clearly, there is no reason to include 
such a candidate in~$S$ 
because every solution with is will be rejected as not satisfying~(IV).

The most time consuming step of the classification 
was from $3$-local to $4$-local codes, 
because the number of nonequivalent
$3$-local codes is huge, see Table~\ref{t:ires}.
Evaluating the CPU time necessary to solve
such amount of exact covering problems
showed that it was
out of reasonable resources.
However, it could be observed that the most
of $3$-local codes have no continuation,
and another software, integer linear programming (ILP) solver called \emph{BOP} from \emph{Google OR-Tools}~\cite{ortools},
 was used to quickly decide 
whether a $3$-local code has a $4$-continuation.
Since it did not provide more than one solution
for each problem, the rest of the work for
continuable $3$-local codes was done 
with \texttt{libexact}. This additional filtering step
took around $1300$ CPU days and speeded up the process 
approximately $10$ times or even better.
Another optimization was related with the partition 
of the classified objects into two subclasses, see Section~\ref{s:sq-sqf}.

\subsection{Isomorph rejection}
\label{s:isom}

An important step, traditionally called the 
\emph{isomorph rejection}, is choosing
non\-equiv\-alent representatives from the set 
of all found solutions, intermediate
or final. 
A standard technique to deal with equivalence, 
see~\cite[Sect. 3.3]{KO:alg}, is to
represent codes by graphs
and use a software \cite{nauty2014} 
to deal with graph isomorphism.
With this technique,
for each code~$A$ one can construct
the \emph{canonically labeled graph} $\Canon(A)$
such that two codes are equivalent
if and only if their canonically 
labeled graphs coincide.
The same software finds the order
of the automorphism group
of the graph, i.e., in our case, 
the order of the automorphism group
of the corresponding code.
This can be used for validating 
the results of the computation, as discussed in Section~\ref{s:validation}.

\subsection
[From (3,3)- to 3-local codes]
{From $(3,3)$- to $3$-local codes}
\label{s:33to3}

While at the first steps of the classification,
it was computationally hard to jump 
from $2$-local to $3$-local codes,
and it required to consider $(2,3)$-local 
codes as an intermediate step, it was found
empirically that for reconstructing codewords
of weight more that $3$ such additional subdivision
is not reasonable anymore. So, after classifying 
$(3,3)$-local codes, 
they were used to classify $3$-local codes
(which differ from $(3,3)$-local codes 
only by the concept of equivalence),
and then work with $r$-local codes, 
$r=3,4,5,\ldots$.

Each equivalence class of $3$-local codes
can be represented at most twice 
in the collection of nonequivalent 
$(3,3)$-local codes. 
The procedure to choose only one representative
for each equivalence class is as follows.
Assume we have a $(3,3)$-local code~$A$.
It has exactly two codewords of weight~$1$,
with~$1$ in the first (condition~(II'))
and some $i$th positions, respectively.
By a coordinate permutation that swaps 
the two positions corresponding to the weight-$1$
codewords (we denote this operation by~$\sigma$),
we get another local code~$\sigma(A)$.
As a $(3,3)$-local code, 
it can be equivalent to $A$ or not.
For $A$ and $\sigma(A)$, we construct canonically labeled
graphs $\Canon(A)$ and $\Canon(\sigma(A))$ (see Section~\ref{s:isom}).
Then, we compare these graph with respect 
to the lexicographic ordering 
(which depends on the computer realization,
but we need an arbitrary fixed 
linear order on the set of graphs).
If $\Canon(A)\le \Canon(\sigma(A))$, then we keep
$A$ in the collection of representatives;
if $\Canon(A) > \Canon(\sigma(A))$, we remove it.
In the last case, the equivalence class 
containing $(3,3)$-local code~$\sigma(A)$
must be represented by some~$B$
in our classification of $(3,3)$-local codes,
where $\Canon(\sigma(A))=\Canon(B)$.
It follows that 
$\Canon(B)<\Canon(\sigma(B))=\Canon(A)$
in this case, which means that we will keep 
$B$ in the collection of representatives
when we process it.
As a result, if the equivalence class of 
$3$-local codes was represented twice, by~$A$ and~$B$, as $(3,3)$-local codes, 
we keep only one of $A$, $B$; if it was represented once, by~$A$, then $\Canon(A) = \Canon(\sigma(A))$, and we keep~$A$.
This procedure is not very fast as it requires 
calculation of two canonical labeled graphs
for each representative, but it does not require 
much memory and is well parallelizable. It took
$124$ days of CPU time.

\subsection{Square and square-free partitions}\label{s:sq-sqf}

We say that a code
or a local code is 
\emph{square} (\emph{square-free}) if it
includes (does not include)
a quadruple $\{\vc x$, $\vc y$, $\vc z$, $\vc v\}$
of vertices equivalent to 
$\{\,0...0\,0011$,
$0...0\,0101$,
$0...0\,1010$,
$0...0\,1100\,\}$.
We see that 
$L_{4,4,4}$ and $L_{4,8}$
are square, while 
$L_{5,7}$, $L_{6,6}$, $L_{12}$
are square-free.
Moreover,
\begin{lemma}
For every
square $\{14;2\}$-code
$C$, there is an equivalent code
that is a continuation
of~$L_{4,4,4}$ or~$L_{4,8}$.
\end{lemma}

Based on this, we separately classify
square codes as continuations of~$L_{4,4,4}$
and~$L_{4,8}$ and square-free codes
as continuations of~$L_{5,7}$, $L_{6,6}$,
and~$L_{12}$ (keeping at each step only square-free
representatives of local codes).

\subsection{Validation}
\label{s:validation}
In each step of the classification,
by double-counting the total number of
found continuations with the help of the orbit-stabilizer
theorem (see~\cite[Sect. 10.2]{KO:alg} 
for the general strategy),
we can
partially validate the results of this step.
In our case, the validation is based 
on the following straightforward fact.

\begin{proposition}\label{p:val1}
Let $C_0$ be an $r$-local code,
and let $C_1$ be 
an $r'$-local 
code continuing~$C_0$, $r'\ge r$.
Then the number of 
$r'$-continuations of $C_0$ 
that are equivalent to~$C_1$ equals
$|\Sym(C_0)|/|\Sym(C_1)|$,
where $\Sym(C_i)$ is the set of all 
coordinate permutations that 
stabilize $C_i$ set-wise.
\end{proposition}

Proposition~\ref{p:val1} 
is used to confirm the results of finding
$(r+1)$-continuations of each given representative
$C_0$
of $r$-local partitions. During the computation,
we find all possible $(r+1)$-continuations 
and hence their number.
After the isomorph rejection, we find
nonequivalent representatives of 
$(r+1)$-continuations, 
and with Proposition~\ref{p:val1}
find again the total number of continuations.
In the case of any random or systematic mistake
in the computing, the results of the calculation 
the number of continuations in two different ways
will be different with very high probability. 
Checking this equality is one of standard ways
to confirm the results of classification 
of combinatorial structures.

The same approach was used to check the results of
the continuation of $(r_0,r_1)$-local codes, 
with the only difference that instead of 
$\Sym(C_i)$ we used 
$\Sym_0(C_i)$ consisting of all permutations
from $\Sym(C_i)$ that fix the first coordinate.

A slightly different realization of the same principle 
was used to validate the final classification of
$\{14;2\}$-codes.
For each such code~$C$, we denote by $\Loc_6(C)$
the set of all $6$-local codes that have 
a continuation equivalent to~$C$.
Nonequivalent representatives of $\Loc_6(C)$
can be easily found as follows:
for each translation $c+C$, $c\not\in C$,
remove all codewords of weight more than~$6$;
keep only nonequivalent representatives.
Checking that all equivalence
classes of such $6$-local partitions
(derived from the resulting $\{14;2\}$-codes)
were found during the classification of 
all $6$-local partitions additionally confirms
that no $6$-local partitions were missed.
Indeed, the set of all possible $6$-local partitions
is the union of $\Loc_6(C)$ over all 
$\{14;2\}$-codes $C$.
If, as the result of a wrong classification,
we found a proper subset of this set, missing the rest,
then the probability that the found subset 
is also the union of $\Loc_6(C)$ over some 
$\{14;2\}$-codes~$C$ is very small.

The validation approach from the paragraph above 
was separately applied for square and square-free
(see the definition in Section~\ref{s:sq-sqf})
$\{14;2\}$-codes~$C$.
The union of $\Loc_6(C)$ over square-free codes~$C$
consists of $110$ equivalence classes.
For square codes, we count only nonequivalent
$6$-local partitions
that continue~$L_{4,4,4}$ or~$L_{4,8}$; their number $5099717$ in 
$\bigcup\Loc_6(C)$ is also in agree
with what was found, see {Table~\ref{t:ires}}. 

\section{Results}\label{s:res}

\begin{table}[ht]
$
\def\dy{\mbox{ days}}
\def\dyy{\mbox{ day}}
\begin{array}{r@{\ }||@{\ }c|@{\ }c }
 &L_{4,4,4},\ L_{4,8} & L_{5,7} ,\,L_{6,6} ,\,L_{12} \mbox{ ($^*$: square-free)} \\ \hline
(2,3)\mbox{-local} & 14 + 59 & 33 + 37 + 196  \\  
(3,3)\mbox{-local} & 73762927 + 1586116921 & 1280242055 + 543652569  + 7755763093  \\
 \mbox{$3$-local: \it all,}
& 36904735 + 793121035 & 640150181 + 271854554 + 3877947089  \\
\mbox{\it square-free,} &    & 166208491^* + 71966561^* + 1014622649^*\\ 
\mbox{\it continuable} & 17044 + 78904 & 25^* + 30^* + 679^*  \\
4\mbox{-local} & 4753786 + 29233429 & 9^* + 0^* + 117^*  \\
5\mbox{-local} & 15286921 + 16399650 & 9^* + 0^* + 101^*  \\
6\mbox{-local} & 1688762 + 3410955 & 9^* + 0^* + 101^*  \\
\end{array}
$
\caption{Intermediate results of the computation}
\label{t:ires}
\end{table}

The numbers of equivalence classes 
of local codes, found as intermediate results,
are shown it Table~\ref{t:ires};
finding $(3,3)$-local codes took 
$339$ 
CPU days, $3$-local codes 
$124$ 
days,
continuable $3$-local codes
$1293$
days,
$4$-, $5$-, and $6$-local codes $50$, 
$115$, 
and $112$ 
CPU days, respectively.

\begin{theorem}[computational]\label{th:main14}
There are $30848$ equivalence classes
of orthogonal arrays OA$(2048,14,2,7)$;
eight of them are square-free.
$14960$ of them ($4$~square-free)
are punctured $1$-perfect codes. 
The total number of different
OA$(2048,14,2,7)$
is $541012580165257200$ ($267743838601839600$ punctured $1$-perfect codes).
\end{theorem}

It is well known
that every binary orthogonal array
OA$(N,n,2,t)$ of even strength~$t$
can be obtained by shortening 
OA$(2N,n+1,2,t+1)$, see e.g. \cite[Theorem~6.1]{BGS:96}.
Consequently, we can obtain 
the classification of OA$(1024,13,2,6)$
from OA$(2048,14,2,7)$, 
by shortening each representative 
in each coordinate and comparing 
the results for equivalence.

\begin{theorem}[computational]\label{th:main13}
There are $202917$ equivalence classes
of orthogonal arrays OA$(1024,13,2,6)$.
$100473$ of them are punctured 
shortened $1$-perfect codes.
\end{theorem}

It is known \cite[Remark~2]{Kro:OA13} that if $C$ is 
a shortened orthogonal array attaining~\eqref{eq:Fried},
then 
$\mathrm{even}(C) \cup (\mathrm{odd}(C)+\overline 1)$
and
$\mathrm{odd}(C) \cup (\mathrm{even}(C)+\overline 1)$
are completely regular codes 
(in general they can be nonequivalent),
in our case, with intersection array $\{12,2;2,12\}$.

Finally, it is straightforward that 
since a $\{14;2\}$-code~$C$ is an independent set,
both $\mathrm{even}(C)$ and $\mathrm{odd}(C)$ are
$\{14,12,2;2,12,14\}$-codes 
(again, they can be nonequivalent), 
and we can classify such codes too.

\begin{theorem}[computational]\label{th:main13crc}
In $Q_{13}$, there are $247904$ equivalence classes
of completely regular codes
with intersection array $\{12,2;2,12\}$.
In $Q_{14}$, there are $36137$ equivalence classes
of completely regular codes
with intersection array $\{14,12,2;2,12,14\}$.
\end{theorem}

\section[Properties of OA(2048,14,2,7)]{Properties of OA$(2048,14,2,7)$}
\label{s:prop}
\subsection[Almost-OA(2048,14,2,7+1)]{Almost-OA$(2048,14,2,7+1)$}
\label{s:aOA}

We say that an orthogonal array
is \emph{almost-OA$(N,n,2,t+1)$}
if it is OA$(N,n,2,t)$
and every subgraph of $Q_n$ isomorphic to
$Q_{n-t-1}$ contains 
$N/2^{t+1}$, $N/2^{t+1}-1$, or $N/2^{t+1}+1$
elements of the array.

To formulate the result of this section,
we need to define a $1$-perfect code~$P$,
obtained by shortening one 
very symmetric extended $1$-perfect 
code~$P'$ in~$Q_{16}$
(a binary code
is called an extended $1$-perfect if puncturing
it in any coordinate results in a $1$-perfect code).
The code~$P'$ is defined
as the orbit of the all-zero word under
the automorphism group of order~$2048$
with generators shown in Table~\ref{t:gens}.
\begin{table}[ht]
\begin{center}
 \makebox[0mm][c]{\GeneratorsTable}
\end{center}
 \caption{Generators of an automorphism group that acts regularly on the extended $1$-perfect code~$P'$. Each automorphism is written in the form $[\vc v,\pi]$,
 where $\vc v$ is a translation vector and
 $\pi$ is a coordinate permutation (the coordinates are represented by hexadecimal digits; the permutations are also indicated by arrows);
 the action of $[\vc v,\pi]$ on the vertices of $Q_{16}$ is $\vc x \to \vc v + \pi(\vc x)$.}\label{t:gens}
\end{table}
 This group acts regularly on~$P'$,
 which means that the code is \emph{propelinear},
 in the sense of~\cite{RifPuj:1997}.
 Puncturing $P'$ in any coordinate results in the same, 
 up to equivalence, $1$-perfect code; to be explicit,
 we denote by~$P$ the result of puncturing $P'$
 in the last coordinate.

\begin{theorem}\label{th:aOA}
Up to equivalence, there is exactly one 
almost-OA$(2048,14,2,7{+}1)$; 
it is the $1$-perfect code~$P$
punctured in one of the first $8$ positions.

The code~$P$ itself is not an 
almost-OA$(2048,15,2,7+1)$ (in a $Q_7$-subgraph of $Q_{15}$, it can have $0$, $7$, $8$, $9$, or $16$ codewords); in particular, there are no almost-OA$(2048,15,2,7{+}1)$ arrays \emph{(which can also be seen from column~$\kappa'_7$ of~\cite[Table~XIX]{OPP:15})}.

There is exactly one almost-OA$(1024,14,2,6+1)$,
which is~$P$ shortened in 
one of the last~$7$ positions.

There are six almost-OA$(1024,13,2,6+1)$.
\end{theorem}

\subsection[Distance-2 connected components]{Distance-$2$ connected components}
\label{s:d2}
\newcommand\TIMES{X}
The \emph{minimum-distance} graph 
of a code is a graph with the codewords as the vertices, two codewords being adjacent
if and only if the distance between them
is minimum over all pairs of codewords. 
Some information about the structure 
of connected components 
of the minimum-distance graph of a $\{14;2\}$-code
was analysed.
Note that such a code is punctured $1$-perfect 
if and only if the minimum-distance graph is bipartite
(to make a $1$-perfect code, one can append one bipartite set with symbol~$0$ and the other bipartite set with symbol~$1$).
Connected components correspond to so-called
$i$-components in $1$-perfect codes of length~$15$, which where analysed in details in~\cite{OPP:15}.
Such a component of a $\{14;2\}$-code can have one of the sizes $1\TIMES$, 
$2\TIMES$, 
$4\TIMES$, 
$6\TIMES$, 
$7\TIMES$, 
$8\TIMES$, where $\TIMES = 128$.
The proportion between sizes of the components
in the even-weight subcode can be one of 
8,
1:7,
2:6,
1:1:6,
4:4,
2:2:4,
1:1:2:4,
1:1:1:1:4,
2:2:2:2,
1:1:2:2:2,
1:1:1:1:2:2,
1:1:1:1:1:1:1:1;
in the odd-weight subcode it is always the same.
In all proportions, the components of size at least
$4\TIMES$ can be non-bipartite, with odd-girth
$5$, $7$, or~$9$ for the size $8\TIMES$ and only~$7$ for smaller sizes (again, all characteristics of the even-weight and odd-weight components coincide, with only one exception with
components of size $8\TIMES$, one bipartite and one not, of odd-girth~$7$).
All components have girth~$4$, except 
some components 
of size~$8\TIMES$ with girth~$5$
and
some bipartite components 
of sizes~$8\TIMES$ and~$7\TIMES$ with girth~$6$.
The even-weight and odd-weight subcodes 
can be nonequivalent; in the bipartite case 
this can only happen with components 
of size~$8\TIMES$ and girth~$4$, 
while in the non-bipartite case,
in all cases except components of size~$7\TIMES$
(there are only $2$ such non-bipartite codes)
or of girth~$5$ ($4$ codes). All components are 
antipodal (contain $\vc x$ and $\vc x +\bar 1$ simultaneously), 
which has a theoretical explanation only 
in the bipartite case,
see~\cite[Sect.~VI]{KroPot:multifold}.

\subsection{Derived group divisible designs}
\label{s:GDD}
A~\emph{group divisible design} 
$(k,\lambda)$-GDD of type~$m^l$
is the triple $(S,G,B)$
from a set $S$ (of~\emph{points}),
a partition~$G$ of~$S$ into $l$ subsets (\emph{groups}) of size~$m$, $|S|=lm$, and a collection~$B$ of $k$-subsets (\emph{blocks}) of~$S$
\emph{such that} two different points
from~$S$ either belong to a common group and do not belong 
to a common block or belong 
to exactly~$\lambda$ common blocks and do not belong to a common group.
If $\vc c$ is a codeword
of a $\{14;2\}$-code~$C$,
$G$ and $B$ are the 
sets of supports
of weight-$2$ and weight-$3$ elements in $\vc c+C$,
respectively,
then $(\{1,...,14\},G,B)$
is a $(3,2)$-GDD design of type~$2^7$, called \emph{derived} from~$C$.

\begin{theorem}
There are $103966$ ($84852$) nonisomorphic
$(3,2)$-GDD designs of type~$2^7$
that are derived from arrays OA$(2048,14,2,7)$ 
(respectively, only from 
punctured $1$-perfect codes).

There are $208$ ($170$) arrays OA$(2048,14,2,7)$
having only one derived $(3,2)$-GDD designs, up to isomorphism.
There are $206$ ($170$) $(3,2)$-GDD designs derived from such arrays.
\end{theorem}

\section*{Acknowledgments}
The author thanks Aleksandr Krotov
for invaluable help with software,
Aleksandr Buturlakin and
Patric \"Osterg{\aa}rd for insightful discussions,
and the Supercomputing
Center of the Novosibirsk State University 
for provided computational resources.

The work was funded by the Russian Science Foundation (22-11-00266), \url{https://rscf.ru/project/22-11-00266/}.

\subsection*{Data availability} 
The dataset containing the results 
of the classifications of OA$(2048,14,2,7)$ and OA$(1024,13,2,6)$ 
is available in the IEEE DataPort
repository~\cite{Perfect-related},
files \\
{H(14,2)\_((0,14)(2,12))*.txt}, and
{H(13,2)\_((9,2,2)(12,0,1)(12,1,0))*.txt} in the archive
{equitable\_partitions*.zip}.
 

\begin{thebibliography}{10}

\bibitem{BGS:96}
J.~Bierbrauer, K.~Gopalakrishnan, and D.~R. Stinson.
\newblock Orthogonal arrays, resilient functions, error-correcting codes, and
  linear programming bounds.
\newblock {\em \href{http://epubs.siam.org/journal/sjdmec}{SIAM J. Discrete
  Math.}}, 9(3):424--452, Aug. 1996.
\newblock \DOI{10.1137/S0895480194270950}.

\bibitem{BRZ:CR}
J.~Borges, J.~Rif\`a, and V.~A. Zinoviev.
\newblock On completely regular codes.
\newblock {\em \href{http://link.springer.com/journal/11122}{Probl. Inf.
  Transm.}}, 55(1):1--45, Jan. 2019.
\newblock \DOI{10.1134/S0032946019010010}.

\bibitem{BMS:2017:few}
P.~Boyvalenkov, T.~Marinova, and M.~Stoyanova.
\newblock Nonexistence of a few binary orthogonal arrays.
\newblock {\em
  \href{http://www.sciencedirect.com/science/journal/0166218X}{Discrete Appl.
  Math.}}, 217(2):144--150, Jan. 2017.
\newblock \DOI{10.1016/j.dam.2016.07.023}.

\bibitem{BulRy:2018}
D.~A. Bulutoglu and K.~J. Ryan.
\newblock Integer programming for classifying orthogonal arrays.
\newblock {\em \href{http://ajc.maths.uq.edu.au}{Australas. J. Comb.}},
  70(3):362--385, 2018.

\bibitem{FDF:CorrImmBound}
D.~G. Fon-Der-Flaass.
\newblock A bound on correlation immunity.
\newblock {\em \href{http://semr.math.nsc.ru}{Sib. \`Elektron. Mat. Izv.}},
  4:133--135, 2007.
\newblock Online: \url{http://mi.mathnet.ru/eng/semr149}.

\bibitem{Friedman:92}
J.~Friedman.
\newblock On the bit extraction problem.
\newblock In {\em Foundations of Computer Science, IEEE Annual Symposium on},
  pages 314--319, Los Alamitos, CA, USA, 1992. IEEE Computer Society.
\newblock \DOI{10.1109/SFCS.1992.267760}.

\bibitem{HSS:OA}
A.~S. Hedayat, N.~J.~A. Sloane, and J.~Stufken.
\newblock {\em Orthogonal Arrays. Theory and Applications}.
\newblock Springer Series in Statistics. Springer, New York, NY, 1999.
\newblock \DOI{10.1007/978-1-4612-1478-6}.

\bibitem{KO:alg}
P.~Kaski and P.~R.~J. {\"O}sterg{\aa}rd.
\newblock {\em Classification Algorithms for Codes and Designs}, volume~15 of
  {\em Algorithms Comput. Math.}
\newblock Springer, Berlin, 2006.
\newblock \DOI{10.1007/3-540-28991-7}.

\bibitem{KasPot08}
P.~Kaski and O.~Pottonen.
\newblock libexact user's guide, version 1.0.
\newblock Technical Report 2008-1, Helsinki Institute for Information
  Technology HIIT, 2008.

\bibitem{Kirienko2002}
D.~Kirienko.
\newblock On new infinite family of high order correlation immune unbalanced
  {B}oolean functions.
\newblock In {\em Proceedings 2002 IEEE International Symposium on Information
  Theory, Lausanne, Switzerland, June 30 -- July 5, 2002}, page 465. IEEE,
  2002.
\newblock \DOI{10.1109/ISIT.2002.1023737}.

\bibitem{Knuth:DLX}
D.~E. Knuth.
\newblock Dancing links.
\newblock E-print cs/0011047, arXiv.org, 2000.
\newblock \DOI{10.48550/arXiv.cs/0011047}.

\bibitem{Perfect-related}
D.~Krotov.
\newblock Perfect and related codes.
\newblock IEEE Dataport, 2022--2023.
\newblock \DOI{10.21227/w856-4b70}.

\bibitem{Kro:21066}
D.~S. Krotov.
\newblock Equitable $[[2,10],[6,6]]$-partitions of the $12$-cube.
\newblock E-print 2012.00038, arXiv.org, 2020.
\newblock Available at \url{http://arxiv.org/abs/2012.00038 }.

\bibitem{Kro:OA13}
D.~S. Krotov.
\newblock On the {OA}(1536,13,2,7) and related orthogonal arrays.
\newblock {\em
  \href{http://www.sciencedirect.com/science/journal/0012365X}{Discrete
  Math.}}, 343(2):111659/1--11, 2020.
\newblock \DOI{10.1016/j.disc.2019.111659}.

\bibitem{Kro:Redundancy23}
D.~S. Krotov.
\newblock Completely regular codes as optimal structures.
\newblock In {\em 2023 {XVIII} {I}nternational {S}ymposium ``{P}roblems of
  {R}edundancy in {I}nformation and {C}ontrol {S}ystems'' ({REDUNDANCY})},
  pages 82--87. IEEE, 2023.
\newblock \DOI{10.1109/Redundancy59964.2023.10330185}.

\bibitem{KroPot:multifold}
D.~S. Krotov and V.~N. Potapov.
\newblock On multifold packings of radius-1 balls in {H}amming graphs.
\newblock {\em
  \href{http://ieeexplore.ieee.org/xpl/RecentIssue.jsp?punumber=18}{IEEE Trans.
  Inf. Theory}}, 67(6):3585--3598, June 2021.
\newblock \DOI{10.1109/TIT.2020.3046260}.

\bibitem{KroVor:2020}
D.~S. Krotov and K.~V. Vorob'ev.
\newblock On unbalanced {B}oolean functions with best correlation immunity.
\newblock {\em \href{http://www.combinatorics.org}{Electr. J. Comb.}},
  27(1):\#P1.45(1--24), 2020.
\newblock \DOI{10.37236/8557}.

\bibitem{nauty2014}
B.~D. McKay and A.~Piperno.
\newblock Practical graph isomorphism, {II}.
\newblock {\em J. Symb. Comput.}, 60:94--112, 2014.
\newblock \DOI{10.1016/j.jsc.2013.09.003}.

\bibitem{OstPot:15}
P.~R.~J. {\"O}sterg{\aa}rd and O.~Pottonen.
\newblock The perfect binary one-error-correcting codes of length $15$: Part
  {I}---classification.
\newblock {\em
  \href{http://ieeexplore.ieee.org/xpl/RecentIssue.jsp?punumber=18}{IEEE Trans.
  Inf. Theory}}, 55(10):4657--4660, 2009.
\newblock \DOI{10.1109/TIT.2009.2027525}.

\bibitem{OPP:15}
P.~R.~J. {\"O}sterg{\aa}rd, O.~Pottonen, and K.~T. Phelps.
\newblock The perfect binary one-error-correcting codes of length $15$: Part
  {II}---properties.
\newblock {\em
  \href{http://ieeexplore.ieee.org/xpl/RecentIssue.jsp?punumber=18}{IEEE Trans.
  Inf. Theory}}, 56(6):2571--2582, 2010.
\newblock \DOI{10.1109/TIT.2010.2046197}.

\bibitem{PWLL:2021}
S.~Pang, J.~Wang, D.~K.~J. Lin, and M.-Q. Liu.
\newblock Construction of mixed orthogonal arrays with high strength.
\newblock {\em Ann. Stat.}, 49(5):2870--2884, 2021.
\newblock \DOI{10.1214/21-AOS2063}.

\bibitem{ortools}
L.~Perron and V.~Furnon.
\newblock {OR-Tools}, v9.8.
\newblock \url{https://developers.google.com/optimization/}.

\bibitem{RifPuj:1997}
J.~Rif\`a and J.~Pujol.
\newblock Translation-invariant propelinear codes.
\newblock {\em
  \href{http://ieeexplore.ieee.org/xpl/RecentIssue.jsp?punumber=18}{IEEE Trans.
  Inf. Theory}}, 43(2):590--598, 1997.
\newblock \DOI{10.1109/18.556115}.

\bibitem{SEN:2010:OA}
E.~D. Schoen, P.~T. Eendebak, and M.~V.~M. Nguyen.
\newblock Complete enumeration of pure-level and mixed-level orthogonal arrays.
\newblock {\em
  \href{http://onlinelibrary.wiley.com/journal/10.1002/(ISSN)1520-6610}{J.
  Comb. Des.}}, 18(2):123--140, 2010.
\newblock \DOI{10.1002/jcd.20236}.

\end{thebibliography}

\providecommand\href[2]{#2} \providecommand\url[1]{\href{#1}{#1}}
  \def\DOI#1{{\href{https://doi.org/#1}{https://doi.org/#1}}}\def\DOIURL#1#2{{\href{https://doi.org/#2}{https://doi.org/#1}}}

\end{document}